\documentclass{article}
\usepackage{amsmath}
\usepackage[russian]{babel}
\usepackage[utf8]{inputenc}
\usepackage{graphicx}
\usepackage{soul}
\usepackage[bottom=2.5cm,left=3cm,right=3cm,top=3cm]{geometry}

\begin{document}

\large
\noindent
{\bf Yashina M.V., Tatashev A.G.}

\section*{\LARGE Approximate Computation of Loss Probability for Queueing System with Capacity Sharing Discipline}

{\normalsize  A multi-channel queueing system is considered. The arriving requests differ in their type. Requests of each type arrive according to a Poisson process. The number of channels required for service with the rate equal to~1 depends of the request type. If a request is serviced with the rate equal to~1, then, by definition, the  length of the request equals to the total service time. If at arrival moment, the idle channels is sufficient, then the arriving request is serviced with the rate~1. If, at the arrival moment, there are no idle channel, then the arriving request is lost. If, at arrival moment, there are idle channels but the number of idle channels is not sufficient for servicing with rate~1, then the request begins to be in service with rate equal to the ratio of the number of idle channels to the number of the channels required for service with the rate~1. If a request is serviced with a rate less than~1 and another request leaves the system, then the service rate increases for the request in consideration.   Approximate formula for loss probability has been proposed. The accuracy of approximation is estimated. Approximate values are compared with exact values found from the system of equations for the related  Markov chain stationary state probabilities.}
\vskip 10pt
{\bf Keywords:} multi-channel systems, capacity sharing, loss probability, approximate approaches.
\vskip 10pt

\section*{1. Introduction}

Models of transmission in information networks can be based on multi-channel queueing systems in that the number of channels required to service a request depends on the request type. This number of channels corresponds to the width of the frequency range required for the transmission of messages. In the multi-channel loss system of this kind, an arriving request is accepted if, at the arrival moment, there is the required number of idle channels, and the request is lost otherwise, [1].  In a queueing system with capacity sharing disciplines, there can be different priority classes for requests, [2].  In such disciplines, incoming requests may be denied in order to reserve a certain number of channels for requests of higher priority classes. In [1], [2], the applications of the models based on queueing system of such class in the analysis of the operation of 4G/5G mobile systems is discussed. 

With a sufficiently large number of channels, an exact computation of the loss of a request (the denial of service or redirection to service in another system) becomes impossible in practice due to the large dimension of the problem. Therefore, approaches to an approximate calculation of this probability are being developed [2]--[4]. In [5], a multi-channel loss system with priority disciplines is studied.  

Other kinds of bandwidth sharing discipline belong to the class of processor (device) sharing disciplines. In the case of the egalitarian processor sharing discipline, which was introduced by L. Kleinrock [6], any number of requests can be in the system, and each of requests is serviced at a rate inversely proportional to the number of the requests. In a more general case of processor sharing discipline, simultaneously serviced requests can be serviced at different rates [7]. In the limited processor sharing discipline [8]--[10], the number of simultaneously serviced requests cannot exceed a prescribed value. If at arrival moment,  there is the maximum permissible number of serviced requests in the system, then the arriving request is lost or goes to the queue, depending on the discipline version.

In this paper, a multi-channel queueing system is described. Requests of different types arrive according to a Poisson process. The request type is characterized by the number of channels required to service the request with the rate~1. If, at arrival moment, there are no idle channels, the arriving request is lost, and, if there are idle channels, but their number is less than the number of channels for servicing the request with the rate~1, the arriving request is serviced with a reduced rate.

In Section 2, the considered system is described. In Section 3, a possible practical interpretation of the considered system is described. In Section 4, the ergodicity of the considered system is proved. In Section 5,  an approximate approach is proposed for calculating the loss probability. The accuracy of approximation is estimated in Section~6. The  approximate values are compared with the exact values found from the system of equations for the stationary probabilities of the related Markov chain states.

\section*{2. Description of system}

Suppose requests of $N$ types arrive to an $m$-channel system according to a Poisson process. For the requests of the type~$i,$ the arrival rate equals $\lambda_i,$ $i=1,\dots,N.$ Suppose, for the service of a request of the type~$i,$ $d_i\le m$ channels are required, $i=1,\dots,N.$ Suppose $B_i(x)$ is the distribution of the service time for the type $i$ request; $b_i<\infty$ is the average value of the length, $i=1,2,\dots,N.$ If there is at least one idle channel but the number of idle channels is not less than $d_i,$ and  a request of 
the type~$i$ arrives, then the arriving request occupies $d_i$ channels and is serviced with the rate~1. If there are $1\le k1<d_i$ idle channels, then the arriving request occupies $k$ channels and begins be serviced with the rate~$\frac{k}{d_i},$ $1\le k<d_i,$ $i=1,2,\dots,N.$ If, in the time interval during that a request of the type~$i$ occupies $k<d_i$ channels, the service of another request ends, and $k_2$ channels, which were servicing this request become idle, then the number of channels servicing the type~$i$ request in consideration increases by $\min(k_2,d_i-k_1),$ and the service rate is equal to $\min(\frac{k_1+k_2}{d_i},1),$ $i=1,2,\dots,N.$  
  
Suppose a request of the length $l$ arrive at time $t_0$ and, at time $t,$ the request is serviced with the rate $\sigma(t).$ Then, at time $t_0+u$ such that 
$$\int\limits_{t_0}^{t_0+u}\sigma(t)dt=l,$$
where $l$ is the request length, the request service ends.  

The distribution of the $i$ type request length is $B_i(x),$ $i=1,\dots,N.$

\section*{3. Possible practical interpretation of system}

Let us give a possible practical interpretation of the system. Suppose there are $m$ robots that upload vehicles, which arrive to an unloading point. There are $N$ types of vehicles. The maximum number of robots that service a vehicle of the type $i$ equals $d_i,$ and, if $N$ robots upload this vehicle, then the unloading time distribution is $B_i(x),$ $i=1,\dots,N.$ If, at the type~$i$ request arrival moment, there are $d_i$ idle robots, then all these robots begin to unload the vehicle, and, if there are $1\le k<d_i$ idle robots, then  $k$ robots unload the vehicle,  $i=1,2,\dots,N.$ Robots that have finished to unload another vehicle join the unloading of the vehicle. The unloading rate is proportional to the number of robots participating in the unloading. If there is no idle robot at the moment of the arrival, the vehicle leaves the unloading point and is unloaded in an alternative way.

\section*{4.  Ergodicity of system}

The expectation of request sojourn time in the system is finite and the number of requests that are in the system simultaneously is limited. Hence the system is ergodic [11], 
i.~e., there exists the stationary distribution of the number of requests in the system (there is the system state distribution), and this distribution does not depend on the initial state of the system. 

 \section*{5. Approximate approach to compute loss probability}

Suppose
$$A_i=\lambda_ib_i:\eqno(1)$$ 
$A_i$ is the incoming load that was created by requests of the $i$-th type under the assumption that each request occupies only one channel;
$$A=\sum\limits_{i=1}^N A_i,\eqno(2)$$
$$\overline{d}=\frac{\sum\limits_{i=1}^N A_id_i}{\sum\limits_{i=1}^N A_i};\eqno(3)$$ 
$\overline{d}$ is the weighted average number of the channels needed for the service of the type~$i$ request with the rate~1;
$$v=\frac{n}{\overline{d}};\eqno(4)$$
$c$ is the stationary state probability.

Suppose
$$d_1=\dots=d_N=\overline{d},\eqno(5)$$
and the value $v$ computed according to (1)--(4) is an integer. Obviously, in this case, the request loss probability is the same as in the system M/G/$v$/0 with arriving load equal 
to $A,$ and hence the probability is computed according to the first Erlang formula
$$c=E=\frac{A^v}{v!\left(1+\sum\limits_{k=1}^n\frac{v^k}{k!}\right)}.$$

Let the value $v=\frac{n}{\overline{d}}$ is fractional. The first Erlang formula may be generalized for a fractional value $v$ [12]
$$c=E(A,v)=\frac{A^ve^{-A}}{\int\limits_A^{\infty}y^me^{-y}dy}.\eqno(6)$$ 
The value of the integral in the right hand side denominator in (6) can be found with aid of the incomplete gamma function tables.

{\it Assume that the value computed according to formulas (2)--(4), (6) is an approximate value of the loss probability for the service system under consideration.
}

Note that, if the value $v$ is integer, than the value computed according to (6) is the same as the value computed according M/G/$v$/0 for the relates loss system. Thus, if (5) holds, (or there is only one type of requests), and the value $v$ computed in accordance with (4) is an integer, then the approach provides the exact value of the loss probability. An approximate value does not depend on the request length distribution under the assumption that the average value is prescribed.

\section*{6. Estimation of approximation accuracy}

In Tables~1 -- 3, approximate values are shown for the loss probability $p$ in a system with two types of requests. In the left part of a cell, the exact value is given. This value is found from the equations for the related Markov chain. In the right side of a cell, an approximate value computed with aid of the proposed approach; $m$ is the number of channels; $d$ is the number channels required for type $i$ service with rate~1; $m$ is the number channels; $\lambda_i$ is the $i$ type request process input rate; $b_i$ is the
type $i$ average length, $i=1,2.$ The first type request length is distributed exponentially.  The second type request length is distributed exponentially (Table~1), or the length distribution is the  the order~2 Erlang distribution, i.~e., the length is the sum of two random values distributed exponentially with the same value (Table~2), i.~e., or the length 
is order~2 hyperexponential distribution, i.~e., with probability~$\alpha_1$,  the length is distributed exponentially with the average $1/\mu_1,$ and, with probability 
$\alpha_ 2,$  the length is distributed exponentially with the average $1/\mu_2,$ and $\alpha_2/\alpha_1=\mu_2/\mu_1$  (balanced hyperexponential distribution).  
\vskip 160pt

{\bf Table 1.} Exact and approximate loss probabilities. The length of a request is distributed exponential with parameter depending on the request type. 
\vskip 10pt
\begin{tabular}{|c|c|c|c|c|c|c|c|c|}
\hline
id&$m$&$d_1$&$d_2$&$\lambda_1$&$\lambda_2$&$b_1$&$b_2$ &$p$\\
\hline
1&2&1&2&1&1&1/2&1/4&0.2632|0.2614\\
\hline
2&2&1&2&2&1&1/3&1/3&0.3289|0.3259\\
\hline
3&2&1&2&1&1&1&1/2&0.4444|0.4405\\
\hline
4&3&2&3&1&1&1/2&1/3&0.3265|0.3007\\
\hline
5&2&1&2&1&9&1/10&1/20&0.3187|0.3139\\
\hline
6&5&1&4&9&9&1/6&1/12&0.4498|0.4080\\
\hline
7&5&1&4&9&9&1/3&1/6&0.5813|0.5683\\
\hline
8&3&2&3&1&1/4&1&1&0.3796|0.4030\\
\hline
9&3&2&3&1&1/2&1&1&0.4698|0.5054\\
\hline
10&3&2&3&1&1&1&1&0.6875|0.6075\\
\hline

\end{tabular}

\vskip 20pt
{\bf Table 2.} Exact and approximate loss probabilities. The length of a type 1 request is distributed exponentially. The type 2 request length distribution is the order~2 Erlang distribution. 

\vskip 20pt
\begin{tabular}{|c|c|c|c|c|c|c|c|c|}
\hline
id&$m$&$d_1$&$d_2$&$\lambda_1$&$\lambda_2$&$b_1$&$b_2$ &$p$\\
\hline
1&2&1&2&1&1&1/2&1/4&0.2644|0.2614\\
\hline
2&2&1&2&2&1&1/3&1/3&0.3301|0.3259\\
\hline
3&2&1&2&1&1&1&1/2&0.4458|0.4405\\
\hline
4&3&2&3&1&1&1/2&1/3&0.3284|0.3007\\
\hline
\end{tabular}

\vskip 70pt
{\bf Table 3.} Exact and approximate loss probabilities. The length of a type 1 request is distributed exponentially. The type 2 request length distribution is the order~2 balanced distribution. 
\vskip 10pt
\begin{tabular}{|c|c|c|c|c|c|c|c|c|}
\hline
id&$m$&$d_1$&$d_2$&$\lambda_1$&$\lambda_2$&$b_1$&$b_2$ &$p$\\
\hline
1&2&1&2&1&1&1/2&1/4&0.2628|0.2614\\
\hline
2&2&1&2&2&1&1/3&1/3&0.3287|0.3259\\
\hline
3&2&1&2&1&1&1&1/2&0.4441|0.4405\\
\hline
4&3&2&3&1&1&1/2&1/3&0.3262|0.3007\\
\hline
\end{tabular}
\vskip 20pt

In Table 4,  the exact and approximate values for the probability $p$ in the system with one type of requests; $m=3$ is the number of channels; $d$ is the number of channels required for the request service with the rate~1; $\lambda$ is the request input rate; $b$ is the average length of request.   
\vskip 20pt
{\bf Table 4.} Exact and approximate values of the loss probability.
{\normalsize
\vskip 10pt
\begin{tabular}{|c|c|c|c|c|c|c|}
\hline
id&Request length distribution&m&$d$&$\lambda$&$b$&$p$\\
\hline
1&exponential&3&2&1&1&0.2500|0.3259\\
\hline
2&order 2 Erlang distribution&3&2&1&1&0.2458|0.3259\\
\hline
3&oder 2 balanced hyperexponential distrbution&3&1&1&1&0.2514|0.3259\\
\hline
4&exponential&19&10&1&1&0.2083|0.2216\\
\hline
5&exponential&9&5&1&1&0.2174|0.2449\\
\hline
6&exponential&7&2&2&1&0.1185|0.1444\\
\hline
7&exponential&11&2&2&1&0.1265|0.1545\\
\hline
\end{tabular}
}

\section*{7. Conclusion}

Approximate values are proposed to estimate loss probability for a multi-channel system. For the service of a request with the rate~1, a number of channels is required, and, if the number of channels is not sufficient, then the request is serviced with a reduced rate.  If there is no idle channel, then the arriving request is lost. An approximate formula is proposed for the loss probability. The accuracy of the approximation is estimated. The values computed in accordance with the approximate formula are compared with the exact  
values found from the system of equations for the stationary probabilities of the related Markov chain.

\section*{References}

\hskip 16pt 1. Stasiak M., Glabowski M., Wisniewski A., Zwierzykowski P. Modeling and dimensioning of mobile wireless networks: from GSM to LTE.  John Wiley \& Sons, 2010.

2. Moscholios I.D., Vassilakis V.G., Logothetis, M.D. Boucouvalas A.C. State-dependent bandwidth sharing policies for wireless multirate loss networks. IEEE Transactions on Wireless Communications, 2017, vol.~16, no.~8, pp.~5481--5497.\\
DOI: 10.1109/TWC.2017.2712153

3. Tatashev A.G., Yashina M.V., Buslaev D.A., Serebryakov E. Evaluation of V2X communication architecture characteristics based on loss queueing system with bandwidth sharing discipline //2022 Systems of Signals Generating and Processing in the Field of on Board Communications, IEEE, 2022,  pp. 1--5.\\
DOI: 10.1109/IEEECONF53456.2022.9744315

4. Yashina M.V., Evgrafov V.V., Tatashev A.G., Yashin V.B. Queueing System with Priority Bandwidth Sharing Discipline for V2X Technology. In 2024 Systems of Signal Synchronization, Generating and Processing in Telecommunications (SYNCHROINFO) (2024)  pp. 1--5. IEEE.
DOI: 10.1109/SYNCHROINFO61835.2024.10617457
 
5. Yashina M., Tatashev A., de Alencar M.S. (2023). Loss probability in priority limited processing queueing system. Mathematical Methods in the Applied Sciences, 46(12), 13279--13288. DOI: 10.1002/mma.9249
 
6. Klienrock L. Time-shared systems: a theoretical treatment // J.~Assoc. Comput. Mach., 1967, vol.~14, no.~2, pp.~242--261. DOI:~10.1145/321386.321388

7. Yashkov S.F. Mathematical issues of the theory of queuing systems with processor sharing // Results of Science and Technology. Series "Probability Theory. Mathematical Statistics. Theoretical Cybernetics", 1990, vol.~29, pp. 3--82. (In Russian.)

8. Telek M., van Houtd B. Response time distribution of a class of limited processor sharing queues // Proceedings of IFiP WG 7.3 Performance Conferences, November 13--17, 2017.

9. Dudin A.N., Dudin S.A., Dudina O.S., Samuylov K.E. Analysis of queueing model with processor sharing discipline and customer impatience // Operation Research Perspectives, vol.~5, pp. 245--255, 2018. DOI: 10.1016/j.orp.2018.08.003

10. Alencar M., Yashina M., Tatashev A. (2021, October). Loss queueing systems with limited processor sharing and applications to communication networks. In 2021 International Conference on Engineering Management of Communication and Technology (EMCTECH) (pp. 1-5). IEEE. DOI: 10.1109/EMCTECH53459.2021.9618978

11. Queueing Systems, Volume I. Wiley, 1974.

12. Shneps M.A. Information distribution systems. Computation methods.  Moscow, Svyaz, 1979. (In Russian.)

\end{document}